\documentclass[a4paper,fleqn]{cas-sc}
\usepackage{etoolbox}
\usepackage{expl3}
\usepackage{capt-of}
\usepackage{amsthm}
\usepackage{amsmath}
\usepackage{amssymb}
\newtheorem{remark}{Remark}
\usepackage{graphicx}
\usepackage{float}
\usepackage[normalem]{ulem}
\usepackage[numbers]{natbib}

\def\tsc#1{\csdef{#1}{\textsc{\lowercase{#1}}\xspace}}
\tsc{WGM}
\tsc{QE}
\tsc{EP}
\tsc{PMS}
\tsc{BEC}
\tsc{DE}

\begin{document}
\let\WriteBookmarks\relax
\def\floatpagepagefraction{1}
\def\textpagefraction{.001}
\shortauthors{Shukai Du et~al.}
\shorttitle{Second scaling regime in high-dimensional approximation}

\title [mode = title]{Second scaling regime in high-dimensional function approximation and curse of dimensionality}                      
\author[1]{Shukai Du}
\author[2]{Thu Le}
\author[2]{Samuel N. Stechmann}

\affiliation[1]{
oraganizationsep={},
    organization={Department of Mathematics, Syracuse University}
}

\affiliation[2]{
oraganizationsep={},
    organization={Department of Mathematics, University of Wisconsin--Madison}
}

\begin{abstract}
High-dimensional function approximation is often subject to the curse of dimensionality, commonly described by a hypercube scaling regime with a large expense of $R^d$ grid points, where $R$ denotes the side length of a hypercube in $d$ dimensions.
Here we show that a second scaling regime, which we call the lattice ball scaling regime, also arises and is relevant for practical applications. In this regime, the expense scales as $d^{R^2}/(R^2!)$ and therefore grows polynomially with $d$. 
This regime arises without additional special assumptions of special types of function spaces or error norms.
The lattice ball scaling regime occurs when $R^2\ll d$, due to the discrete nature of counting the expense as the number of lattice points inside a hypersphere. Here we investigate this regime, which arises from the classical problem of approximating high-dimensional functions from common spaces of analytic or $m$-times continuously differentiable functions using orthogonal series (spectral method) such as Fourier, Chebyshev, and Legendre polynomial series.  
For many applications, the expense in the lattice ball scaling regime can be substantially smaller than that of the hypercube scaling regime, making accurate approximations more accessible.
Numerical illustrations include classifying images of handwritten digits from the MNIST dataset and learning the solution operators of partial differential equations, such as Burgers’ equation and the radiative transfer equation.
The potential for connections with artificial neural networks and other machine learning methods is also discussed.
\end{abstract}



\maketitle
\noindent

\section{Introduction}
\subsection{Setting}
\label{subsec:setting}
We consider the approximation of high-dimensional functions using
common spectral bases, such as Fourier, Chebyshev, and Legendre polynomial
bases. We first introduce the general setting by considering a target
function $f:\Omega\to\mathbb{R}$, where
$\Omega\subset\mathbb{R}^d$. Standard choices for $\Omega$
include the periodic domain $[-\pi,\pi]^d$ and the hypercube
$[-1,1]^d$.

Let $\{\phi_{\mathbf{k}}\}_{\mathbf{k}\in\mathcal{I}}$ be an
orthonormal tensor-product spectral basis for $L^2(\Omega,\mu)$,
where $\mu$ is the associated product measure and
$\mathcal{I}\subseteq\mathbb{Z}^d$ is the index set, with each
multi-index $\mathbf{k}=(k_1,\ldots,k_d)$
labeling a spectral mode. We write
$
\phi_{\mathbf{k}}(\mathbf{x})
=
\prod_{i=1}^d \varphi_{k_i}^{(i)}(x_i),
$ 
where $\mathbf{x}=(x_1,\ldots,x_d)$ and
$\varphi_{k_i}^{(i)}$ denotes the one-dimensional basis function
indexed by $k_i$ in the $i$-th coordinate.

For $f \in L^2(\Omega,\mu)$, its spectral expansion is given by
\begin{equation}
f(\mathbf{x}) = \sum_{\mathbf{k}\in\mathcal{I}} \hat{f}_{\mathbf{k}} \phi_{\mathbf{k}}(\mathbf{x}),
\qquad
\hat{f}_{\mathbf{k}} = \int_{\Omega} f(\mathbf{x}) \overline{\phi_{\mathbf{k}}(\mathbf{x})} \, d\mu(\mathbf{x}).
\label{eq:spectral-expansion}
\end{equation}

\subsection{Hypercube scaling and curse of dimension}

For high-dimensional approximations, 
the cost is commonly described as 
growing exponentially with dimension $d$, as
\begin{equation}
    \operatorname{Vol}([-R,R]^d)=(2R)^d.
    \label{eq:hypercube-volume}
\end{equation}
This is the volume of a hypercube with 
side length $2R$ in $d$ dimensions. 
This cost arises either from a multi-dimensional numerical mesh
with $2R$ grid points in each dimension,
or from a spectral (orthogonal series) expansion
as in \eqref{eq:spectral-expansion}
where a similar multi-dimensional numerical mesh
is used for the multi-dimensional wavevector $\mathbf{k}$.

This situation is often described as
the curse of dimensionality,
since the cost increases rapidly (exponentially)
as the dimension $d$ increases.

\subsection{Hyperbolic cross and other special assumptions}

Many interesting studies have identified ways of 
overcoming the curse of dimension in certain situations.
For instance, some interesting work has identified
specialized settings 
which explicitly or implicitly involve an aspect that is sparse or low-rank
(e.g., \cite{wang2025operator})
or which are amenable to the application of artificial neural networks
which show reduced dependence on dimensionality
(e.g., \cite{chen2023deep}).

Other interesting studies have introduced special assumptions such as
special types of function spaces or special versions of error norms,
such as Barron space (e.g., \cite{barron1993universal,e2022barron,lu2021apriori})
or the hyperbolic cross of the sparse grid method or Smolyak method
(e.g., \cite{bungartz2004sparse,dung2018hyperbolic,adcock2022sparse}).

Here, in what follows, we will consider the setting in
section~\ref{subsec:setting} 
without making additional assumptions about special types of function spaces or error norms.

\subsection{Hypersphere scaling: promise and apparent paradox}

While the hypercube and hyperbolic cross are commonly used
in describing the cost of high-dimensional function
approximation, it may be appropriate and natural to
consider the volume of a hyperball, 
i.e., the volume of the interior of the hypersphere.
The hypersphere arises from theoretical decay 
estimates for the magnitude of the coefficients 
$\hat{f}_\mathbf{k}$ in \eqref{eq:spectral-expansion}. 
In such estimates (e.g., \cite{canuto2006spectral,grafakos2014classical,trefethen2017multivariate}), 
the magnitude $|\hat{f}_\mathbf{k}|$ decays 
according to the Euclidean norm $\|\mathbf{k}\|_2$.
Consequently, it is natural, based on these
theoretical considerations, to consider all
wavevectors $\mathbf{k}$ such that
$\|\mathbf{k}\|_2\le R$,
which suggests that the cost should scale as
the volume of $B_R^d$, the ball of radius $R$ in
$d$ dimensions:
\begin{equation}
\operatorname{Vol}(B_R^d) = \frac{\pi^{d/2}}{\Gamma(d/2 + 1)} R^d.
\label{eq:hypersphere-volume}
\end{equation}

At first glance, it seems that the hyperball volume
in \eqref{eq:hypersphere-volume} should be promising and should offer
a massive reduction in cost, since it is well-known
that the hyperball volume will shrink rapidly relative to
the volume in \eqref{eq:hypercube-volume} of its bounding hypercube \cite{vershynin2018high,wainwright2019high}.
For instance, setting $d = 1000$ and $R = 5$, the  volume of the hypercube $[-R,R]^d$ is $10^{1000}$, whereas the volume inside the hypersphere is  
$
\operatorname{Vol}(B_5^{1000}) \approx 2.88 \times 10^{-187}
$,
which is many orders of magnitude smaller.

However,
this immense discrepancy creates an apparent paradox: if the volume inside the hypersphere ($\approx 10^{-187}$) is less than one, does this imply that fewer than a single wavevector $\mathbf{k}$ is needed for the orthogonal series in \eqref{eq:spectral-expansion}? Clearly, something is amiss.

\section{Second scaling regime: the lattice scaling regime}

\begin{figure}[pos=H]
    \centering
    \includegraphics[width=0.25\linewidth]{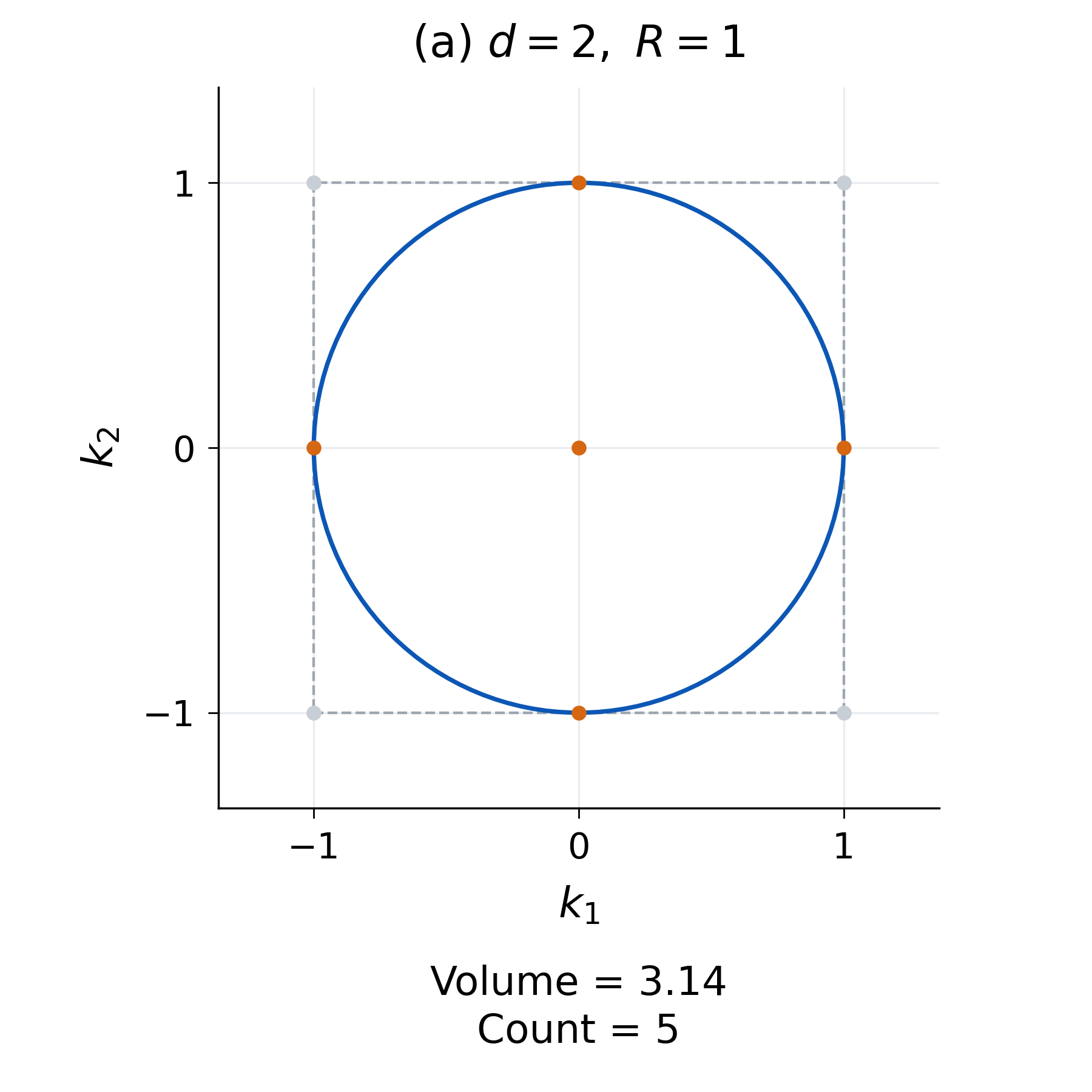}
    \hspace{0.5cm}
    \includegraphics[width=0.25\linewidth]{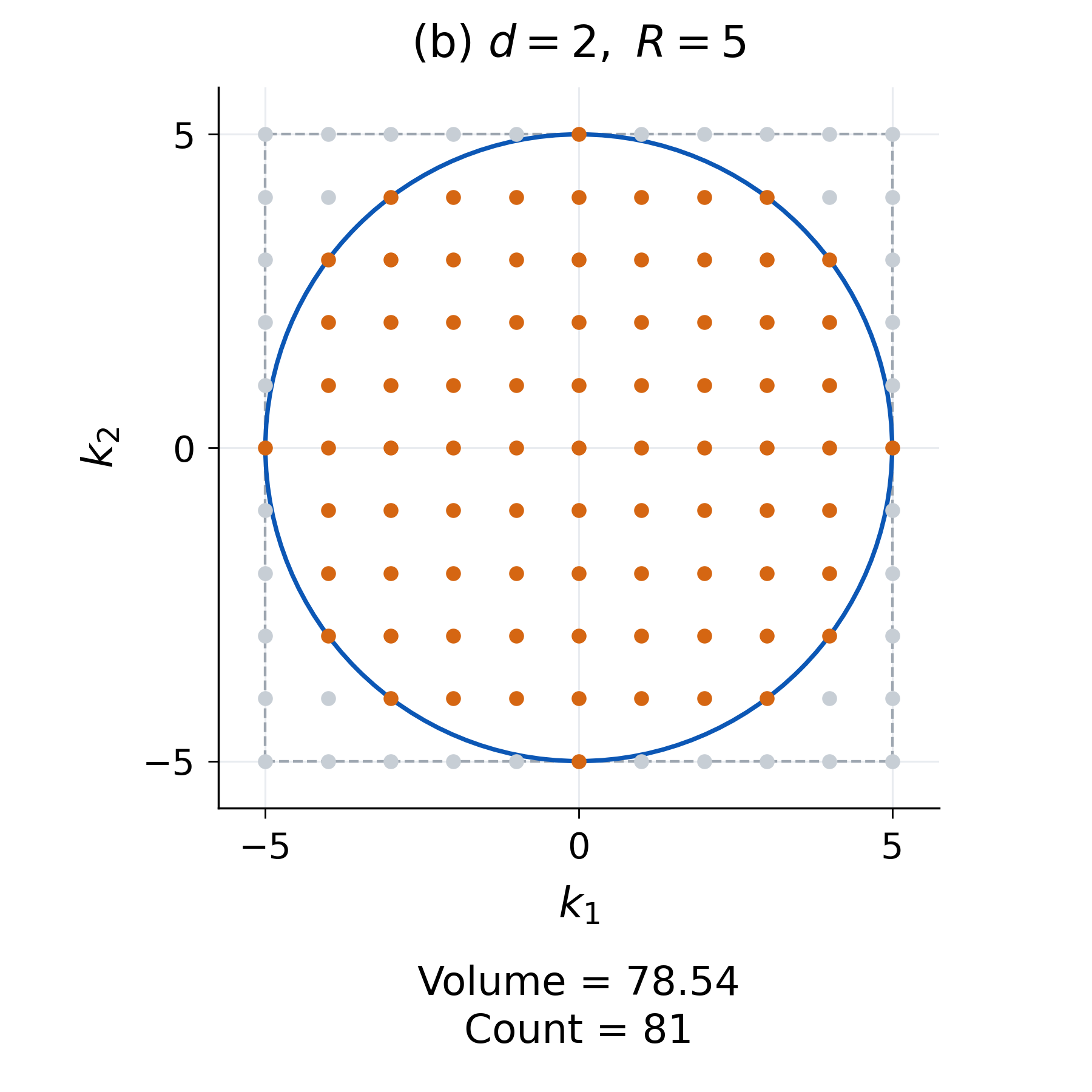}
    \hspace{1cm}
    \includegraphics[width=0.25\linewidth]{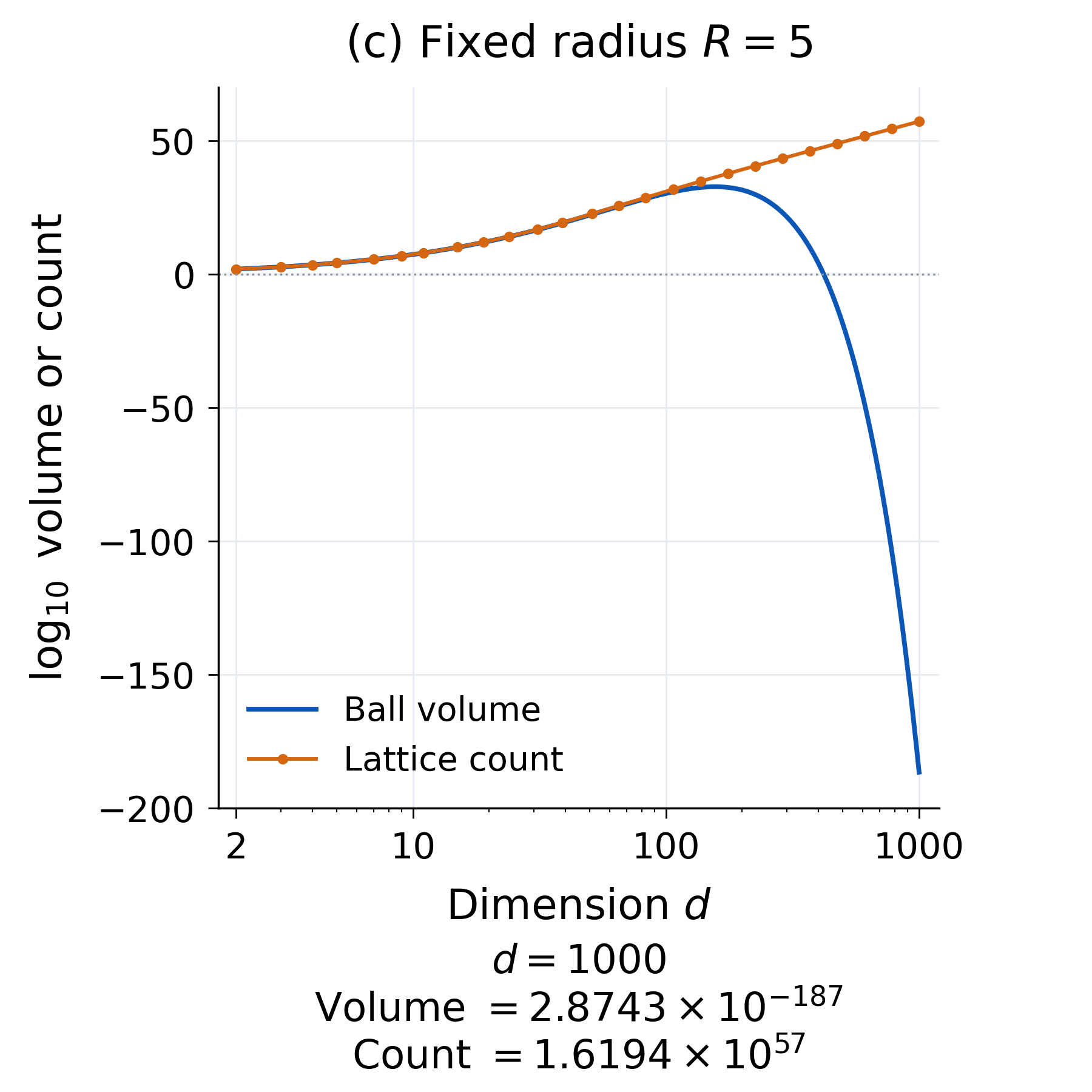}
    \caption{
    Continuous hypersphere volume versus lattice-point cardinality.
    (a)--(b) Two-dimensional balls with $R=1$ and $R=5$;
    orange points indicate lattice points.
    (c) Volume and exact lattice count versus dimension for $R=5$.
    }
    \vspace{-0.5cm}
    \label{fig:ball-volume-lattice}
\end{figure}

In the apparent paradox above, the fallacy lies in tacitly
assuming that hypersphere volume always serves as a direct
proxy for the computational cost.
Since the  wavevectors $\mathbf{k}$ have integer
coordinates, counting them amounts to counting lattice points
inside the hypersphere.
Figure~\ref{fig:ball-volume-lattice} illustrates the distinction
between  volume and this discrete count.

To explain this discrepancy, we distinguish two asymptotic regimes:
\begin{enumerate}
    \item {Large-radius regime (fixed $d$, $R\to\infty$).}
    When the radius $R$ is large relative to the dimension $d$, continuous volume \eqref{eq:hypersphere-volume} provides an asymptotically accurate
estimate of the number of lattice
points in the ball. 

    \item {Lattice scaling regime (fixed $R\ge1$,
    $d\to\infty$).}
     When the dimension $d$ is much larger than the radius $R$,  continuous volume  \eqref{eq:hypersphere-volume}   tends to zero, while the number of lattice
points in the ball continues to grow. Thus, volume no longer provides an accurate estimate
of the discrete lattice-point count.
\end{enumerate}

To quantify the cardinality growth in the lattice scaling regime, we next define finite downward-closed index sets (e.g., \cite{adcock2022sparse,trefethen2017multivariate,hecht2026multivariate}). A set $\Lambda \subset \mathbb{N}_0^d$ is downward-closed if $\mathbf{k} \in \Lambda$ and $0 \le j_i \le k_i$ for all $i = 1, \ldots, d$ imply $\mathbf{j} \in \Lambda$.
We consider downward-closed index sets of the form
\begin{align}
\Lambda_{R,p} = \left\{ \mathbf{k} \in \mathbb{N}_0^d: \|\mathbf{k}\|_p \le R \right\},
\qquad R \ge 0, \quad 1 \le p \le \infty.
\label{eq:setup-index-cutoff}
\end{align}
The common choices are the total-degree ($p=1$),
Euclidean-degree ($p=2$), and maximum-degree ($p=\infty$)
index sets.
Geometrically, these sets consist of the nonnegative lattice
points in an $\ell^1$ ball, an $\ell^2$ ball (enclosed by a
hypersphere), and a hypercube, respectively.
The associated approximation space and orthogonal projection
are defined by
\begin{equation}
V_{\Lambda_{R,p}} = \operatorname{span} \{ \phi_{\mathbf{k}} : \mathbf{k} \in \Lambda_{R,p} \},
\qquad
P_{\Lambda_{R,p}} f = \sum_{\mathbf{k} \in \Lambda_{R,p}} \hat{f}_{\mathbf{k}} \phi_{\mathbf{k}}.
\end{equation}
The orthogonal projection $P_{\Lambda_{R,p}}f$ is the best
approximation to $f$ in $V_{\Lambda_{R,p}}$ with respect to
the $L^2(\mu)$ norm:
\begin{equation}
\|f - P_{\Lambda_{R,p}} f\|_{L^2(\mu)}^2 = \min_{\{c_{\mathbf{k}}\}}  \| f - \sum_{\mathbf{k} \in \Lambda_{R,p}} c_{\mathbf{k}} \phi_{\mathbf{k}} \|_{L^2(\mu)}^2.
\label{eq:setup-tail}
\end{equation}

The maximum-degree and total-degree
index sets have the exact cardinalities:
\begin{equation}
|\Lambda_{R,\infty}|=(R+1)^d,
\qquad
|\Lambda_{R,1}|=\binom{d+R}{R}.
\label{eq:sizelmax}
\end{equation}

For the Euclidean index set, we obtain bounds that determine
the leading-order growth of its cardinality.  Since every
$\mathbf{k}\in\Lambda_{R,2}$ satisfies
$
\|\mathbf{k}\|_1
\le \|\mathbf{k}\|_2^2
\le R^2.
$
Therefore,
$\Lambda_{R,2}\subseteq\Lambda_{R^2,1}$, giving the upper bound
$
|\Lambda_{R,2}|
\le \binom{d+R^2}{R^2}.
$

For $d\ge R^2$, a lower bound follows by considering indices
with exactly $R^2$ coordinates equal to $1$ and all remaining
coordinates equal to $0$. There are $\binom{d}{R^2}$ such
admissible indices.
Consequently,
\begin{equation}
\binom{d}{R^2}
\le |\Lambda_{R,2}|
\le \binom{d+R^2}{R^2}.
\label{eq:l2esitmate}
\end{equation}
For fixed cutoff $R$, both bounds have the same leading term as
$d\to\infty$, yielding
\begin{equation}
|\Lambda_{R,2}|
=\frac{d^{R^2}}{(R^2)!}
 +\mathcal{O}\left(d^{R^2-1}\right)
\sim \frac{d^{R^2}}{(R^2)!}.
\label{eq:regime2}
\end{equation}
Thus,  the Euclidean index set cardinality
grows polynomially in $d$, whereas the
maximum-degree index set cardinality grows exponentially
in $d$.

\begin{remark}
We examine approximation errors for common classes of target
functions and relate them to the truncation radius $R$ and the
coefficient cost associated with different index sets $\Lambda_{R,p}$. We first
consider a function  in the Sobolev class $f\in H^s$ with $s>0$,
for which the following error bounds hold:
\begin{equation}
\|f-P_{\Lambda_{R,p}}f\|_{L^2}
\le
\|f\|_{H^s}
\begin{cases}
(1+R^2/d)^{-s/2}, & p=1,\\[2pt]
(1+R^2)^{-s/2},   & p=2,\infty.
\end{cases}
\label{eq:remark-error-radius}
\end{equation}
Indeed,
$
\|f-P_{\Lambda_{R,p}}f\|_{L^2}^2 =
\sum_{\|\mathbf{k}\|_p>R}
(1+\|\mathbf{k}\|_2^2)^{-s}
(1+\|\mathbf{k}\|_2^2)^s
|\widehat f_{\mathbf{k}}|^2
\le
(1+r_p^2)^{-s}\|f\|^2_{H^s}.
$
Since $\|\mathbf k\|_1\le\sqrt d\,\|\mathbf k\|_2$ and $\|\mathbf k\|_\infty\le\|\mathbf k\|_2$, it follows that $\|\mathbf k\|_p>R$ implies $\|\mathbf k\|_2>r_p$, where $r_1=R/\sqrt d$ and $r_2=r_\infty=R$.

Euclidean-degree and maximum-degree truncations satisfy
the same error bound, but their coefficient costs $C_p:=|\Lambda_{R,p}|$ can differ
substantially.
For fixed $R$,
\eqref{eq:sizelmax} and \eqref{eq:regime2} give
\[
C_2
\sim
\frac{d^{R^2}}{(R^2)!},
\qquad
C_\infty
=
(R+1)^d,
\qquad \text{ as } d\to\infty.
\]
See the supplementary materials for the $p=1$ case. In what follows we focus on $p=\infty$ and $p=2$.

We can also express the approximation error bound \eqref{eq:remark-error-radius} in terms of the cost. 
For maximum-degree truncation, substituting $R=C_\infty^{1/d}-1$ to \eqref{eq:remark-error-radius} gives
\begin{equation}
\|f-P_{\Lambda_{R,\infty}}f\|_{L^2}
\le
\|f\|_{H^s}
(1+(C_\infty^{1/d}-1)^2)^{-s/2} \lesssim C_\infty^{-s/d}.
\label{eq:remark-tensor-cost}
\end{equation}
For Euclidean truncation, the counting estimates and Stirling’s formula give
$
\log C_2
\sim R^2\log\left(\frac{d}{R^2}\right)$ in the lattice scaling regime. 
Asymptotic inversion yields
$
R^2
\sim
\frac{\log C_2}{\log(1+d/\log C_2)}.
$
Substituting into \eqref{eq:remark-error-radius}, we obtain
\begin{equation}
\|f-P_{\Lambda_{R,2}}f\|_{L^2}
\le
(1+o(1))\|f\|_{H^s}
\left(
1+\frac{\log C_2}{\log(1+d/\log C_2)}
\right)^{-s/2}\lesssim
\left({\log C_2^{1/\log d}}\right)^{-s/2}.
\label{eq:remark-lattice-cost}
\end{equation}

For finitely smooth target functions
$f\in C^m(\overline{\Omega})$, where $m\ge1$ is an integer,
the bound in \eqref{eq:remark-error-radius} applies with $s=m$.
 For target functions analytic in a complex neighborhood of
$\overline{\Omega}$ (e.g., \cite{hecht2026multivariate,trefethen2017multivariate}), a similar argument yields,  
for some $c>0$,
\begin{equation*}
\|f-P_{\Lambda_{R,p}}f\|_{L^2}
\lesssim
\begin{cases}
\exp\bigl(-c C_\infty^{1/d}\bigr),
& p=\infty,\\[5pt]
\exp\bigl(-c\sqrt{\tfrac{\log C_2}
{\log(1+d/\log C_2)}}\bigr),
& p=2.
\end{cases}
\label{eq:remark-analytic-cost}
\end{equation*}
\end{remark}

\begin{remark}
   The index cutoff \eqref{eq:setup-index-cutoff} treats all
coordinates equally. Anisotropy can be incorporated through
weighted index sets
\begin{equation}
\Lambda^{\mathbf w}_{R,p}
=
\left\{
\mathbf k\in\mathbb N_0^d:
\|W\mathbf k\|_p\le R
\right\}, \quad 
W=\operatorname{diag}(w_1,\ldots,w_d),
\quad w_j>0,\quad j=1,\ldots,d, \quad 1\le p\le\infty.
\label{eq:weighted_index}
\end{equation}
These sets remain downward closed.  
The weights could chosen in a variety of ways,
such as $w_j=(1+j)^{r}$ with algebraic growth rate $r>0$.
With such weighted spaces, it is the rate $r$
that has an important influence on cost and error,
instead of or in addition to a finite dimension $d$.
See supplementary materials for further discussion,
and see numerical experments below in section~\ref{sec:numerical}
for numerical examples with anisotropic, weighted spaces.

\end{remark}

\section{Numerical examples}
\label{sec:numerical}
We investigate high-dimensional function approximation through three numerical examples: a radial Runge function,
a viscous Burgers solution map, and MNIST classification. 
A fourth example, the radiative transfer equation, is included in supplementary material.
\subsection{Runge function}
\label{subsec:runge}
We first consider the radial Runge function 
$
f(\mathbf x)
 =(1+\alpha\sum_{j=1}^d x_j^2)^{-1},
$ with $\mathbf x\in[-1,1]^d$ and $\alpha=1$. 
All coordinates enter symmetrically, providing an isotropic
benchmark for the index cutoff  $\Lambda_{R,p}$
in \eqref{eq:setup-index-cutoff}.

We approximate the Runge function using the tensor-product Legendre expansion described in Section~2.
We draw sample inputs from the
uniform distribution on $[-1,1]^d$ and estimate the
expansion coefficients by solving the empirical
least-squares of \eqref{eq:setup-tail}
on the training samples. The approximation error
is evaluated on the test samples.  
 
\begin{figure}[pos=H] 
    \centering
    \includegraphics[width=0.9\linewidth]{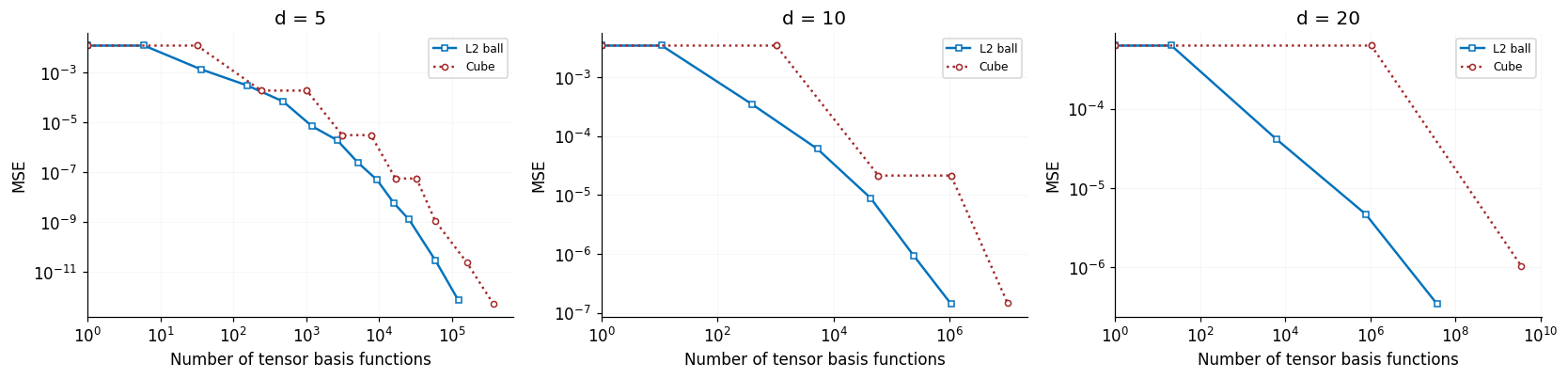}
    \caption{Approximation error versus cost for
    the radial Runge function
    $f(\mathbf x)=(1+\sum_{j=1}^d x_j^2)^{-1}$
    on $[-1,1]^d$, with $d=5,\,10,\,20$.
  Hypersphere truncation achieves comparable accuracy
with substantially fewer basis functions than hypercube
 truncation.}
    \label{fig:runge-errors}
    \vspace{-0.5cm}
\end{figure}

In the lattice scaling regime,  the number of basis functions retained by
hypersphere   truncation ($p=2$) grows polynomially with $d$, whereas
hypercube truncation  ($p=\infty$) exhibits exponential growth.
The error-versus-cost comparison in
Figure~\ref{fig:runge-errors} examines this reduction in
cost and its impact on approximation accuracy.
The results show that hypersphere truncation achieves comparable accuracy with substantially fewer basis functions than hypercube truncation.

\subsection{Burgers equation}
We next consider the viscous Burgers equation
$u_t+uu_x=\nu u_{xx}$ on $[0,2\pi]$, with viscosity $\nu>0$, periodic boundary
conditions,  and initial condition
$u(x,0)=u_0(x)$.  For a fixed observation point \(x_\ast\) and time \(T>0\), we approximate the solution operator \(F: u_0\mapsto u(x_\ast,T)\).

We generate random initial conditions of the form
$u_0(x)=\sum_{i=1}^{K}
[a_i\cos(ix)+b_i\sin(ix)] 
$. 
Each sample pairs the input coefficient vector
$\mathbf x=(a_1,b_1,\ldots,a_{K},b_{K})\in\mathbb R^d$,
where $d=2K$, with the corresponding  target
$y=u(x_\ast,T)$, computed analytically using the
Cole--Hopf transformation.
After rescaling the inputs to $[-1,1]^d$, we approximate  the solution map by solving the empirical least-squares problem \eqref{eq:setup-tail} using tensor-product Legendre polynomials.   See further details of the experimental setup in the supplementary material.  

Figure~\ref{fig:burgers-weighted} compares approximation
error versus  cost with $\nu=0.08$,
$x_\ast=\pi$, and $T=0.5$.
At $d=10$, with $W=I$ in \eqref{eq:weighted_index},  hypersphere truncation 
achieves comparable accuracy with fewer basis functions
than hypercube truncation.

At $d=20$, we also consider weighted hypersphere truncation.
We assign larger weights to modes with smaller amplitudes,
$w_{2i-1}=w_{2i}=(B_{\max}+\varepsilon)/(B_i+\varepsilon)$,
where
$B_{\max}=\max_{1\le i\le K,\;s\in\mathcal I_{\mathrm{train}}}
\max\{|a_i^{(s)}|,|b_i^{(s)}|\}$, and $\varepsilon>0$.
In this example, the input amplitudes decay quickly with
mode number, so the weighted index set retains fewer basis terms
involving the higher frequency modes. Figure~\ref{fig:burgers-weighted}(b)
shows that the weighted truncation further reduces the error at
comparable cost, illustrating the benefit of adapting the index set
to input anisotropy.

\begin{figure}[pos=H]
    \centering
    \begin{tabular}{@{}c@{\hspace{0.02\linewidth}}c@{}}
        \includegraphics[width=0.3\linewidth]{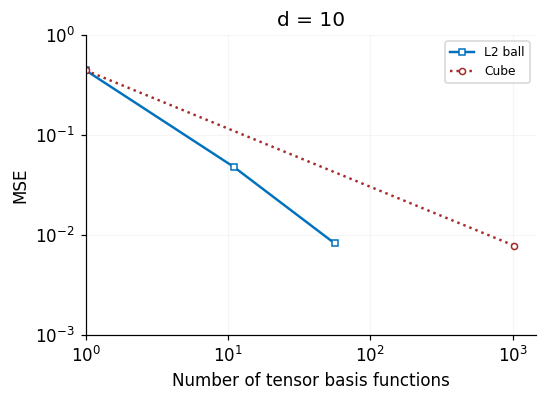}
        &
        \hspace{1cm}
        \includegraphics[width=0.3\linewidth]{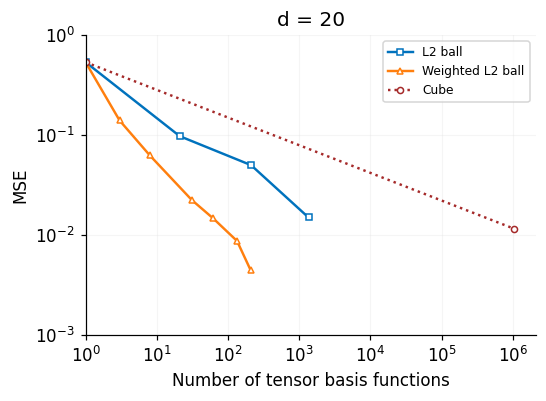}
        \\
        \small (a) & \small (b)
    \end{tabular}
    \caption{Approximation error versus cost for approximation of the viscous
 Burgers solution map $u_0\mapsto u(x_\ast,T)$, with
    $x_\ast=\pi$, $T=0.5$, and $\nu=0.08$.
    (a) (a) At $d=10$, unweighted hypersphere truncation achieves
comparable accuracy with substantially fewer basis
functions than hypercube truncation.
(b) At $d=20$, weighted hypersphere further truncation reduces
the error relative to its unweighted counterpart
at comparable cost.}
    \label{fig:burgers-weighted}
    \vspace{-0.5cm}
\end{figure}
\subsection{MNIST data}
We now consider the classification of handwritten digits in the MNIST dataset,
where each image has $28\times 28$ pixels for a dimension of 784.
An image is represented by a vector
$\mathbf v\in[0,1]^{784}$ with an associated label
$\ell\in\{0,\ldots,9\}$. From the training images, we compute the leading $d$
principal directions and normalize the resulting PCA
coordinates to obtain the inputs $\mathbf z$. Each label is represented by
$\mathbf e_\ell\in\mathbb R^{10}$, with $1$ in the
corresponding component and $0$ elsewhere.  We approximate the target function 
$G(\mathbf z)=\mathbb E[\mathbf e_L\mid\mathbf Z=\mathbf z]$ using tensor-product Legendre polynomials on the training samples.
Here,  
$G_\ell(\mathbf z)=\Pr(L=\ell\mid\mathbf Z=\mathbf z)$.
The corresponding  classifier is
$\ell_\ast(\mathbf z)
=\operatorname*{arg\,max}_{\ell }
G_\ell(\mathbf z)$. 
Further details of the experimental setup are provided in the
supplementary material.

Figure~\ref{fig:MNIST}(a) shows the approximation error per class entry  for the two hypersphere truncations with $d=20$. Both ball truncations achieve reasonably small
approximation errors.  The improvement from weighting is relatively modest in this
example, possibly reflecting a slower decay in input coordinates or the possibility that low-variance components may still carry useful information for classification. We omit hypercube truncation because fitting the resulting basis is computationally expensive.
Figure~\ref{fig:MNIST}(b) shows predictions from the  models on  nine randomly selected test images. Both hypersphere models correctly classify all nine test samples.
\begin{figure}[pos=H]
    \centering
    \begin{tabular}{@{}c@{\hspace{0.02\linewidth}}c@{}}
        \includegraphics[width=0.3\linewidth]{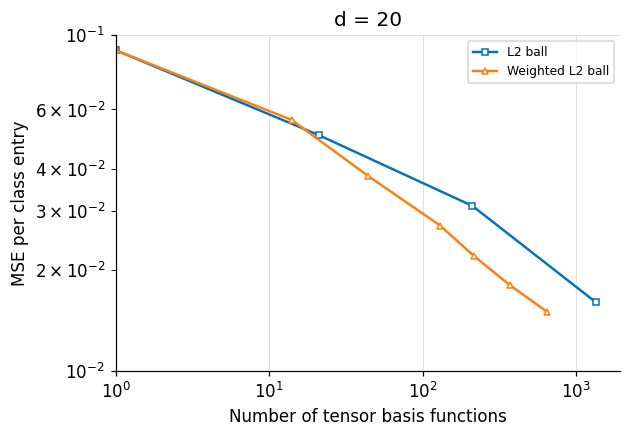}
        &
        \hspace{1cm}
        \includegraphics[width=0.28\linewidth]{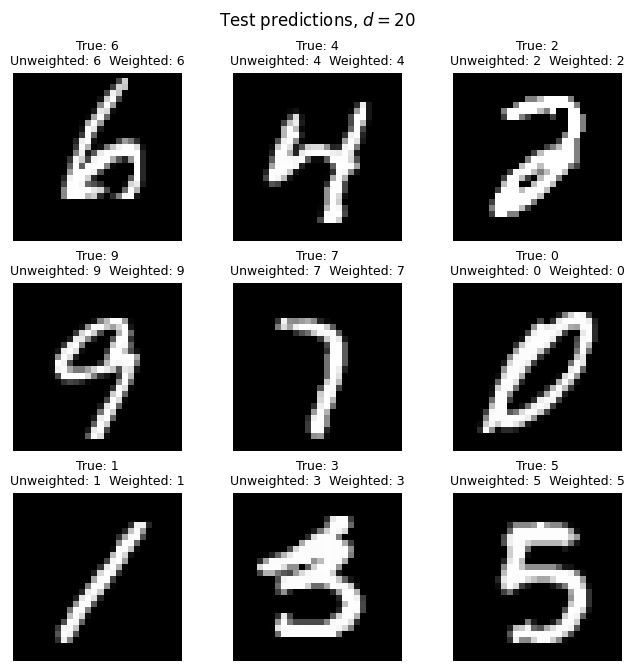}
        \\
        \small (a) & \small (b)
    \end{tabular}
    \caption{MNIST classification with $d=20$ principal components.
    (a) MSE per class entry versus cost for unweighted and weighted hypersphere truncation.
    (b) Predictions on nine randomly selected test images, both hypersphere models correctly classify all nine images.}
    \label{fig:MNIST}
    \vspace{-0.5cm}
\end{figure}
\section{Concluding discussion}
In summary, a second scaling regime was investigated 
in the setting of high-dimensional function approximation.
This lattice ball scaling regime arises for truncation of
orthogonal series with truncation radius $R$ in high dimensions $d$,
when $R^2\ll d$.
In this regime,
the scaling for the cost in \eqref{eq:regime2} is $d^{R^2}/(R^2)!$,
and the error scaling in \eqref{eq:remark-lattice-cost} is more favorable than the curse of dimension
that arises from hypercube scaling.
Several numerical experiments suggest that the cost is feasible
for applications.

For future work, it would be interesting to connect and compare the viewpoint herein
with machine learning approaches such as neural networks.
For example, it would be interesting to investigate
whether the second scaling regime arises in neural networks with training
via stochastic gradient descent.
Results on these topics will be
presented elsewhere in the near future.

\bibliographystyle{cas-model2-names}
\bibliography{cas-refs}

\newpage
\section*{Supplementary Material}

This supplementary material extends the approximation results in
Remarks~1 and~2 of the main text. It also presents an additional
numerical example involving the radiative transfer equation and
provides further details on the experiments reported in Section~3
of the main text.

\section{Additional approximation estimates}
\subsection{Error versus cost for total-degree truncation}
For total-degree truncation ($p=1$), the coefficient cost is
$C_1=|\Lambda_{R,1}|=\binom{d+R}{R}$ for integer $R\ge 1$.
In the lattice scaling regime $R=o(d)$, Stirling's formula gives
\begin{align*}
\log C_1
=R\log d-\log(R!)+O(R^2/d)
\sim R\log\left(\frac{d}{R}\right).
\end{align*}
Writing $L=\log C_1$, we have
$d/L\sim(d/R)/\log(d/R)$. Since $d/R\to\infty$,
$\log(1+d/L)\sim\log(d/R)$, and hence
$
R\sim\frac{\log C_1}{\log(1+d/\log C_1)}.
$
Substituting this expression into the $p=1$ bound in Remark\~1 of the main text yields the following error-versus-cost estimate:
\begin{align}
  \|f-P_{\Lambda_{R,1}}f\|_{L^2}
\le (1+o(1))\|f\|_{H^s}
\left(
1+\frac{(\log C_1)^2}{d[\log(1+d/\log C_1)]^2}
\right)^{-s/2}.  
\end{align}

\subsection{Weighted approximation for anisotropic inputs}
The weighted index sets $\Lambda^{\mathbf w}_{R,p}$ introduced in
Remark~2 of the main text allow us to
account for different rates of coefficient decay across coordinates.
In what follows, we choose the weights for analytic functions on
shrinking input intervals using the radii of the corresponding
Bernstein ellipses. The resulting error bounds decay exponentially in
the truncation radius $R$, with constants independent of $d$ when the
interval widths shrink sufficiently fast.

Let $f=f(x_1,\ldots,x_d)$ be defined on
$U_a=\prod_{i=1}^d[-a_i,a_i]$, where $a_i>0$.
With $x_i=a_i y_i$ and $y_i\in[-1,1]$, set
$$
g(\mathbf y):=f(a_1y_1,\ldots,a_dy_d).
$$
Assume that $g$ extends holomorphically to the Bernstein polyellipse
$E_{\rho_1}\times\cdots\times E_{\rho_d}$, where $\rho_i>1$, and that
$$
\sup_{\mathbf z\in E_{\rho_1}\times\cdots\times E_{\rho_d}}
|g(\mathbf z)|\le M,
$$
where $M$ is independent of $d$.

We use the tensor-product Chebyshev expansion considered in
\cite{hecht2026multivariate,trefethen2017multivariate}.
The same truncation argument applies to other tensor-product polynomial
bases, such as Legendre, given the corresponding coefficient decay estimates.
$$
g(\mathbf y)=\sum_{\mathbf k\in\mathbb N_0^d}
c_{\mathbf k}T_{\mathbf k}(\mathbf y),
\qquad
T_{\mathbf k}(\mathbf y)=\prod_{i=1}^d T_{k_i}(y_i).
$$
Let $m(\mathbf k)=|\{i:k_i>0\}|$. By Cauchy's formula after the
substitution $y_i=(z_i+z_i^{-1})/2$, we have
$$
|c_{\mathbf k}|
\le M\,2^{m(\mathbf k)}\prod_{i=1}^d\rho_i^{-k_i}
= M\,2^{m(\mathbf k)}e^{-\mathbf w\cdot\mathbf k},
\qquad
w_i=\log\rho_i.
$$

Suppose that the input satisfying $a_i\asymp(i+1)^{-r}$ for some $r>0$ and that
$f$ has a bounded holomorphic extension to a product of complex
neighborhoods of the intervals $[-a_i,a_i]$, each of a common
width $\delta>0$. For any $0<\delta_0<\delta$, we can choose
$\rho_i=(\delta_0+\sqrt{a_i^2+\delta_0^2})/a_i$.
Indeed, the scaled ellipse $a_iE_{\rho_i}$ has semiminor axis
$\delta_0$ and lies within distance $\delta_0$ of $[-a_i,a_i]$.
It follows that

$$
\rho_i\asymp(i+1)^r,
$$

As a result,
$
w_i=\log\rho_i\sim\log\bigl((i+1)^r\bigr)
 \text{ as }i\to\infty.$ Consider the weighted index set for $p=2$ or $p=\infty$, and write
$P_R=P_{\Lambda^{\mathbf w}_{R,p}}$ for the chosen $p$.
Since
$
\|W\mathbf k\|_\infty\le \|W\mathbf k\|_2
\le \|W\mathbf k\|_1=\mathbf w\cdot\mathbf k,
$
every omitted index satisfies $\mathbf w\cdot\mathbf k>R$.
Using $\|T_{\mathbf k}\|_\infty\le 1$, for any $0<\theta<1$
we obtain
$$
\begin{aligned}
\|g-P_Rg\|_{L^\infty}
&\le M e^{-\theta R}
\sum_{\mathbf k\in\mathbb N_0^d}
2^{m(\mathbf k)}
e^{-(1-\theta)\mathbf w\cdot\mathbf k}\\
&=M e^{-\theta R}
\prod_{i=1}^d\left(1+\frac{2}{\rho_i^{1-\theta}-1}\right).
\end{aligned}
$$
Writing the product as $C_{d,\theta}$, we have
$$
\log C_{d,\theta}\lesssim
\sum_{i=1}^d(i+1)^{-q},
\qquad q=r(1-\theta).
$$
If $q>1$, the infinite sum converges to $\zeta(q)-1$, where
$\zeta$ is the Riemann zeta function. Thus, $C_{d,\theta}\le C_\theta = O(exp(\zeta(r(1-\theta)))$, where  $C_\theta >0 $ is independent of $d$. Futhermore,
$$
\|g-P_Rg\|_{L^\infty}
\le M C_\theta e^{-\theta R}.
$$
When $r(1-\theta)>1$, or equivalently, $0<\theta<1-\frac1r$.
This range is nonempty when $r>1$.

Similarly,  we have 
$$
\|g-P_Rg\|_{L^2(\mu)}
\le M C'_\theta e^{-\theta R},
\qquad 0<\theta<1-\frac{1}{2r},
$$
where $C'_\theta$ is independent of $d$. This range is nonempty
when $r>1/2$.

\section{Additional numerical results and experimental details}

 \subsection{Runge function}

\begin{figure}[pos=H]
    \centering
    \includegraphics[width=0.9\linewidth]
    {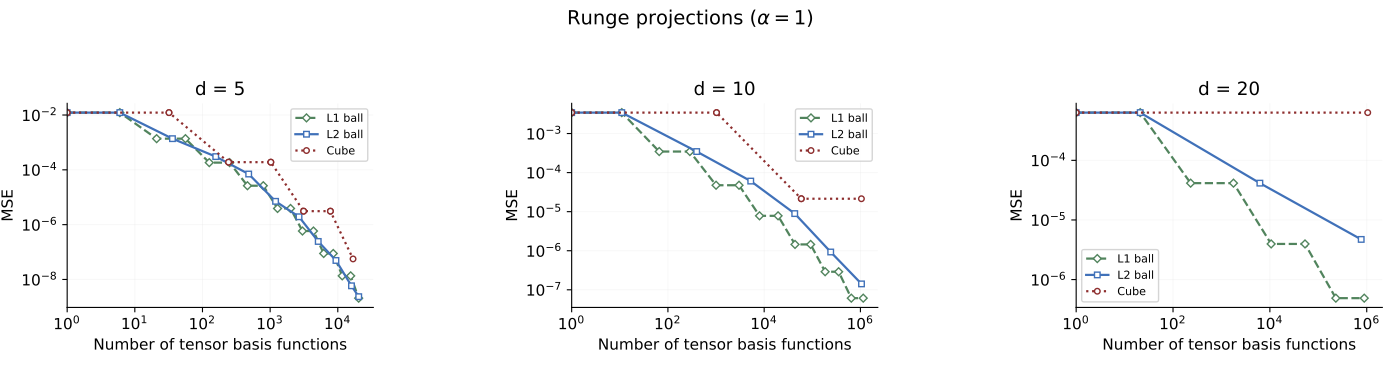}
    \caption{Error versus cost for the radial Runge function
    with $\alpha=1$ and $d=5,10,20$.
    }
    \label{fig:runge_alpha1}
    \vspace{-0.5cm}
\end{figure}

\begin{figure}[pos=H]
    \centering
    \includegraphics[width=0.9\linewidth]
    {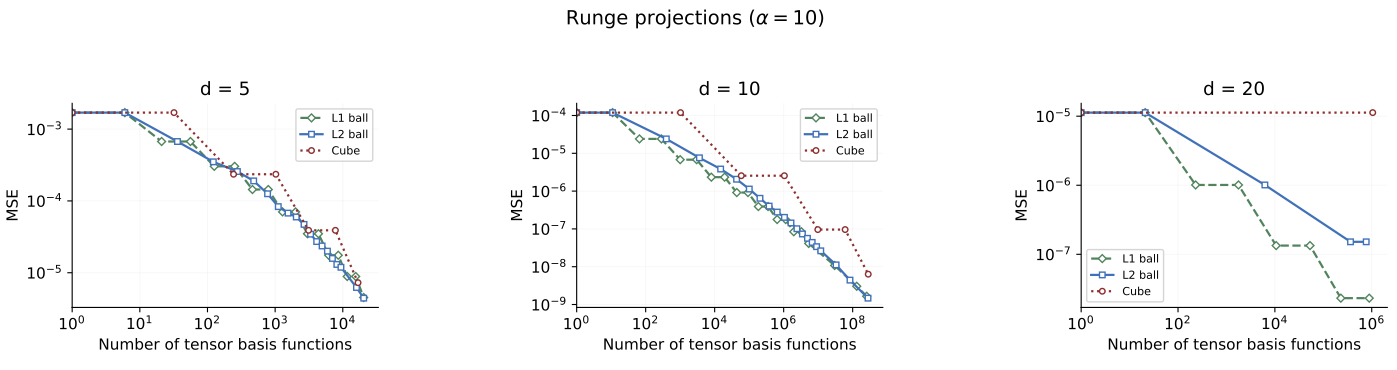}
    \caption{Error versus cost  for the radial Runge function
    with $\alpha=10$ and $d=5,10,20$.}
    \label{fig:runge_alpha10}
    \vspace{-0.5cm}
\end{figure}

We further examine the radial Runge function
$
f(\mathbf{x})
=\left(1+\alpha\sum_{j=1}^{d}x_j^2\right)^{-1},\, \mathbf{x}\in[-1,1]^d,
$
using the tensor-product Legendre projection onto
the index sets $\Lambda_{R,p}$ for $p=1,2,\infty$. Figures~\ref{fig:runge_alpha1} and~\ref{fig:runge_alpha10}
compare the projection errors against the number of basis
functions for the three truncations. The $\ell^1$ and $\ell^2$ truncations generally achieve
smaller errors than the hypercube truncation at comparable cost,
particularly in higher dimensions. 
Although the $\ell^1$ truncation has a weaker error bound
in terms of the radius $R$ than the $\ell^2$ and $\ell^\infty$
truncations,
it achieves smaller errors at comparable cost for $\alpha=1$. For $\alpha=10$, its
performance is comparable to or better than that of the
$\ell^2$ truncation in this setting.

\subsection{Burgers equation}
\begin{figure}[pos=H]
    \centering
    \begin{tabular}{@{}c@{\hspace{0.02\linewidth}}c@{}}       \includegraphics[width=0.4\linewidth]{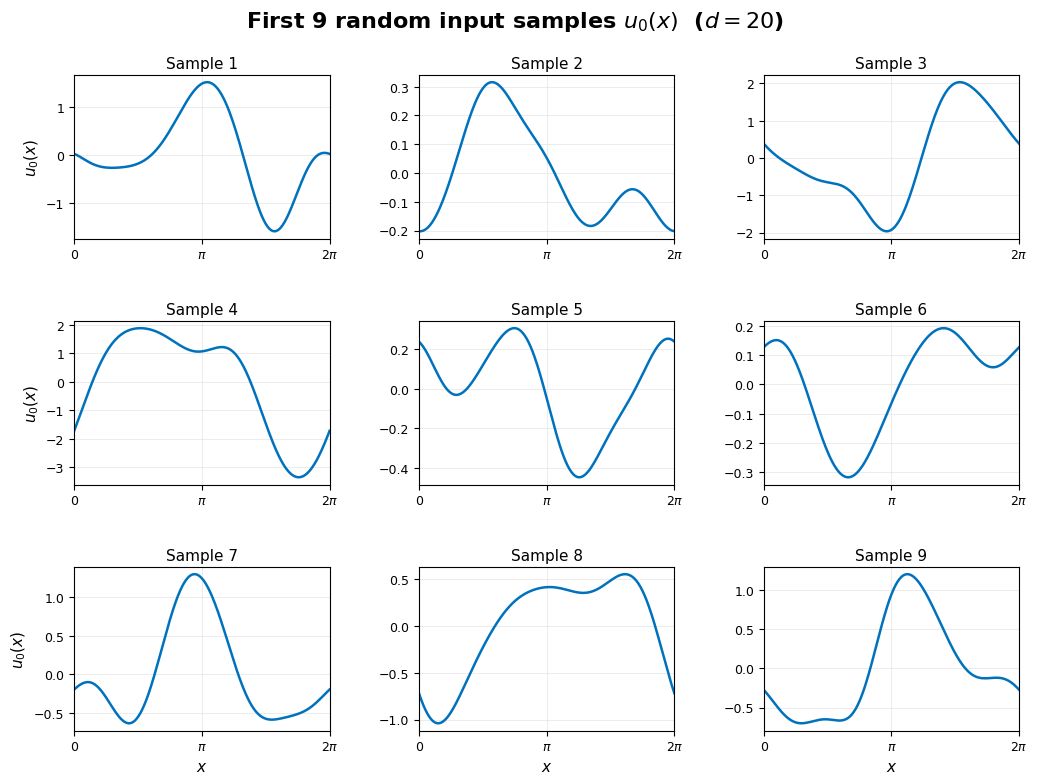}
        &
\hspace{0.3cm}        \includegraphics[width=0.4\linewidth]{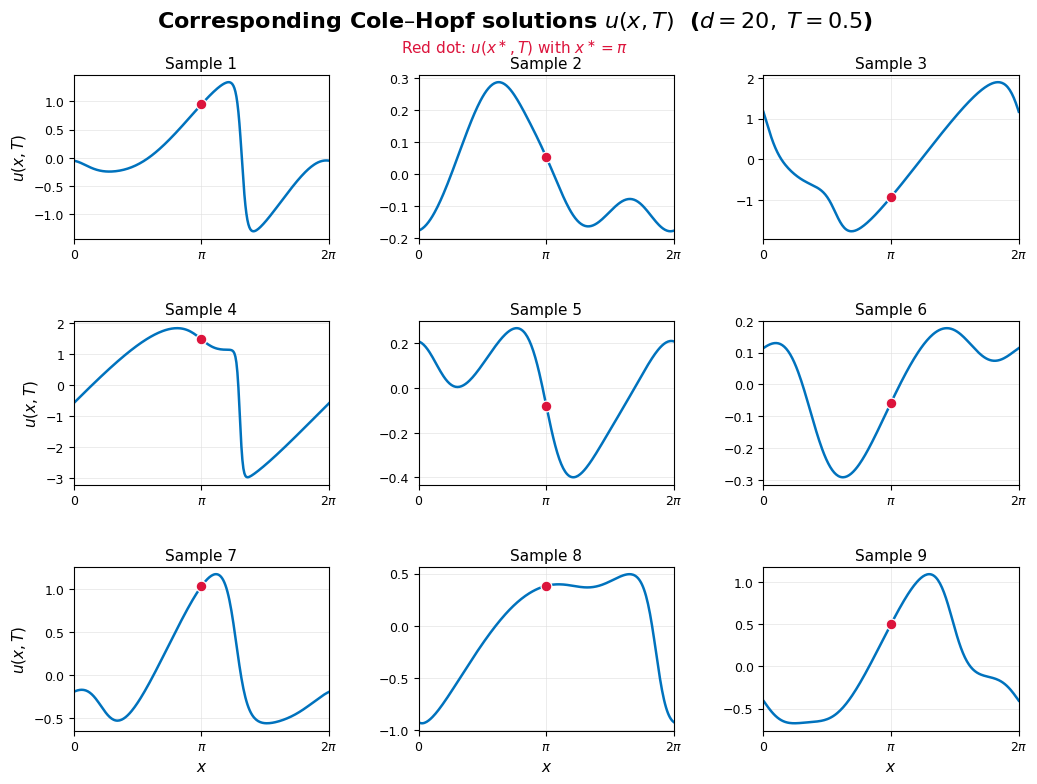}
        \\
        \small (a) & \small (b)
    \end{tabular}
    \caption{Nine randomly generated samples for the viscous Burgers equation
    with input dimension $d=20$.
    (a) Randomly generated initial conditions $u_0(x)$.
    (b) Corresponding solutions $u(x,T)$ at $T=0.5$, computed using the
    Cole--Hopf transformation. The red dots indicate the scalar
    outputs $u(x^*,T)$ at $x^*=\pi$.}
    \label{fig:burgers_samples}
    \vspace{-0.5cm}
\end{figure}

As described in Section~3.2 of the main text, we approximate
the viscous Burgers solution map from its initial condition
to the solution value.
Random initial conditions are generated as
\[
u_0(x)=\sum_{i=1}^{K}
[a_i\cos(ix)+b_i\sin(ix)],
\qquad
a_i=A\rho_i\xi_i,\quad b_i=A\rho_i\eta_i.
\]
For each sample, the amplitude $A$ is drawn uniformly from
$[A_{\min},A_{\max}]$, independently of the random variables
$\xi_i,\eta_i$, which are sampled uniformly from $[-1,1]$.
We use $\rho_i=\widetilde{\rho}_i/\widetilde{\rho}_1$, where
$\widetilde{\rho}_i=\exp[-c_{\mathrm{sm}}(i/K)^{p_{\mathrm{decay}}}]$;
$c_{\mathrm{sm}}$ and $p_{\mathrm{decay}}$ control the decay
of higher modes. See Figure \ref{fig:burgers_samples}(a).

Each sample pairs the input coefficient vector
$\mathbf x=(a_1,b_1,\ldots,a_{K},b_{K})\in\mathbb R^d$,
where $d=2K$, with the  target
$y=u(x_\ast,T;\mathbf x)$, computed using a periodic
Cole--Hopf solver (see Figure \ref{fig:burgers_samples}(b)).
 We normalize the  inputs as
$\mathbf{z}=\mathbf{x}/M_{\mathrm{global}}\in[-1,1]^d$,
where $M_{\mathrm{global}}=A_{\max}\max_i\rho_i$.
We then approximate the solution
map using a tensor-product Legendre expansion:
\[
P_{\Lambda^\mathbf{w}_{R,p}}F  
=
\sum_{\mathbf k\in\Lambda^\mathbf{w}_{R,p}}
c_{\mathbf k}\,\phi_{\mathbf k}(\mathbf z),
\]
with coefficients estimated by solving the empirical
least-squares  on the
training samples. We generate $5000$ sample pairs, using $4000$ for training and $1000$ for testing.

For weighted truncation, we estimate the amplitude of physical
mode $i$ from the training samples and assign the same weight
to the corresponding components of the multi-index:
\begin{align}
w_i
 =\frac{B_{\max}+\varepsilon}{B_i+\varepsilon}, \quad
B_i
 =\max_{s\in\mathcal I_{\mathrm{train}}}
\max\{|a_i^{(s)}|,|b_i^{(s)}|\},
\qquad B_{\max}=\max_i B_i.
\end{align}

The experiments use
$
A_{\min}=0.1,A_{\max}=2,\,
c_{\mathrm{sm}}=20,\,
p_{\mathrm{decay}}=2,
\, \varepsilon=10^{-10},\, x_\ast=\pi,\, \, T=0.5
$ and $\nu=0.08$.  

\subsection{MNIST classification}
As described in Section~3.3 of the main text,
each image is represented by a vector
$\mathbf v\in[0,1]^{784}$ with an associated label
$\ell\in\{0,\ldots,9\}$. From the training images, we compute
the mean $\overline{\mathbf v}$ and the leading $d$ principal
directions $\mathbf u_1,\ldots,\mathbf u_d$. 
The PCA coordinates are $x_j=\mathbf u_j^T(\mathbf v-\overline{\mathbf v})$, $j=1,\ldots,d$.
We then normalize the coordinates
$\mathbf z=\mathbf x/M_{\mathrm{global}}\in[-1,1]^d$,
where
$M_{\mathrm{global}}
=\max_{s\in\mathcal I_{\mathrm{train}}}\max_j|x_j^{(s)}|.$ 

For an image labeled $\ell$, the target is the one-hot vector
$\mathbf e_\ell\in\mathbb R^{10}$, with $1$ in the component
corresponding to class $\ell$ and $0$ elsewhere.
Let $L$ and $\mathbf Z$ denote the random label and normalized
PCA input, respectively. Since
$[\mathbf e_L]_\ell=\mathbf1_{\{L=\ell\}}$, the 
regression target is
\[
G(\mathbf z)=\mathbb E[\mathbf e_L\mid\mathbf Z=\mathbf z],
\qquad
G_\ell(\mathbf z)=\Pr(L=\ell\mid\mathbf Z=\mathbf z).
\]
The ideal classifier is given by
$
\ell_\ast(\mathbf z)
=\operatorname*{arg\,max}_{\ell\in\{0,\ldots,9\}}
G_\ell(\mathbf z).
$ Here, $G$ is the vector of conditional class
probabilities and we estimate $G$ using the spectral expansion:
\[
P_{\Lambda^\mathbf{w}_{R,p}} G (\mathbf z)
=
\sum_{\mathbf k\in\Lambda^\mathbf{w}_{R,p}}
\mathbf c_{\mathbf k}\phi_{\mathbf k}(\mathbf z),
\qquad \mathbf c_{\mathbf k}\in\mathbb R^{10},
\]
with coefficients obtained by solving the empirical least
squares on the training samples.  We approximate using a tensor-product
Legendre expansion.
For weighted truncation, the training PCA eigenvalues
determine
the weights
$w_j=(\sqrt{\mu_1}+\varepsilon)/(\sqrt{\mu_j}+\varepsilon)$,
where $\varepsilon>0$.
This heuristic choice penalizes higher polynomial degrees in lower-variance coordinates.

\section{Radiative transfer equation}
We consider the steady radiative transfer equation (RTE):
\[
\begin{aligned}
\mathbf s\cdot\nabla_{\mathbf x}u
+\sigma_e u
-\sigma_s\int_{S^1}p(\mathbf s,\mathbf s')\,
u(\mathbf x,\mathbf s')\,\mathrm d\mathbf s'
&=0 &&\text{in }\Omega\times S^1,\\
u&=g &&\text{on }\Gamma_-.
\end{aligned}
\]
Here $u(\mathbf x,\mathbf s)$ is the radiative intensity at position
$\mathbf x$ in direction $\mathbf s\in S^1$, $\sigma_e$ and $\sigma_s$
are the extinction and scattering coefficients, and $p$ is the
scattering phase function. The inflow boundary is
\[
\Gamma_-=\{(\mathbf x,\mathbf s)\in\partial\Omega\times S^1:
\mathbf s\cdot\mathbf n(\mathbf x)<0\},
\]
where $\mathbf n$ is the outward unit normal and $g$ is the
prescribed inflow intensity.

With a fixed single-scattering albedo $\omega$, the relation
$\sigma_s=\omega\sigma_e$ allows us to parameterize the optical
properties by a single spatially varying coefficient. Its discrete
representation is
\[
\sigma_h(x_1,x_2)
=\sum_{i=0}^{p}\sum_{j=0}^{p}\hat \sigma_{ij}L_i(x_1)L_j(x_2),
\]
where $L_i$ is the Legendre polynomial of degree $i$. We order the
modal coefficients $\hat \sigma_{ij}$ by increasing total spatial degree $i+j$
to form $\hat \sigma \in\mathbb R^d$. 
We take $p=6$, so the input dimension is $d=(p+1)^2=49$.

We also discretize the
angular variable using $N_a=28$ uniform intervals. For details on how the optical-coefficient samples are generated,
see \cite{DuStechmann2026}.
We use a single scale for all input coordinates. For the reordered
vector $\hat \sigma =(\hat \sigma_1,\ldots,\hat \sigma_d)$, we define
\[
B_j=\max_{s\in\mathcal I_{\mathrm{train}}}|\hat \sigma_j^{(s)}|,
\qquad
M_{\mathrm{global}}=\max_{1\le j\le d}B_j,
\qquad
\mathbf x=  \hat \sigma/M_{\mathrm{global}}.
\]
The normalized training inputs therefore lie in $[-1,1]^d$.

For each input $\hat \sigma$, let $G_{\mathrm{D2N}} $ be the
modal matrix of the local discrete Green operator, which maps the prescribed inflow
intensity to the local interior intensity.
We approximate the map $F$ that takes the  modal optical
coefficients  to the component of
$G_{\mathrm{D2N}}$ along a fixed reference pattern $P$. For each sample, we set 
\[
y=F(\mathbf z)
=\operatorname{Re}\langle
G_{\mathrm{D2N}} ,P\rangle_F,
\]
where $\langle\cdot,\cdot\rangle_F$ denotes the Frobenius inner
product and $P=G_{\mathrm{ref}}/\|G_{\mathrm{ref}}\|_F$ is a fixed reference
pattern for all samples, with $G_{\mathrm{ref}}$ be the matrix from the first stored sample.

We generate $1000$ sample pairs and split them into
$800$ training samples and $200$ test samples. We approximate $F$ by a tensor-product Legendre expansion
 and fit the coefficients by   least squares on the training
samples.  We use $w_j=1$ for the unweighted $\ell^1$ truncation and
\[
w_j=\frac{M_{\mathrm{global}}+\varepsilon}{B_j+\varepsilon},
\qquad \varepsilon=10^{-10},
\]
for the weighted truncation. A smaller training amplitude $B_j$
gives a larger weight and thus raises the cost of increasing the
polynomial degree in coordinate $j$.

\begin{figure}[pos=H]
    \centering
    \includegraphics[width=0.3\linewidth]{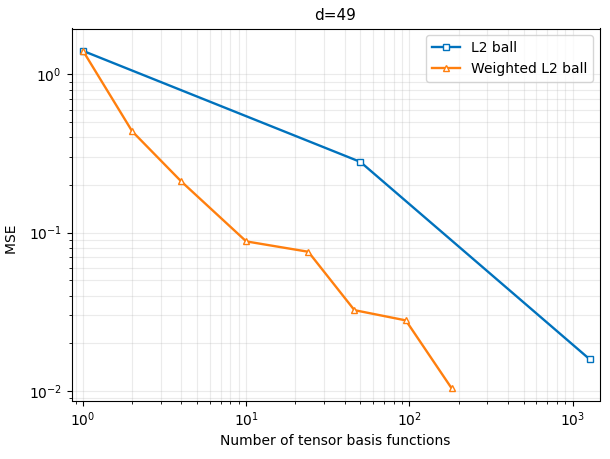}
    \caption{Approximation error versus cost for the RTE map
$F\colon \hat{\sigma}\mapsto
\operatorname{Re}\langle G_{\mathrm{D2N}}(\hat{\sigma}),P\rangle_F$
with $d=49$.  Both  ball truncations attain
small approximation errors. By accounting for input anisotropy, the
weighted truncation achieves comparable accuracy with fewer basis
functions. We omit cube truncation from the figure because it is computationally prohibitive.}
    \label{fig:EL}
\end{figure}

Figure~\ref{fig:EL} shows that both  ball truncations attain
small approximation errors. By accounting for input anisotropy, the
weighted truncation achieves comparable accuracy with fewer basis
functions. We omit hypercube truncation because fitting the resulting
basis is computationally expensive.

\end{document}